\documentclass[12pt,twoside]{article}
\usepackage[english]{babel}
\usepackage{amsmath}
\usepackage{amssymb,amsfonts}
\usepackage{graphicx}
\usepackage[shortlabels]{enumitem}
\usepackage{float}
\usepackage{caption}
\usepackage{titling}
\usepackage{amsthm}

\usepackage{times,amssymb,amscd}
\newtheorem{thm}{Theorem}[section]

\numberwithin{equation}{section}

\newtheoremstyle{mylemma} 
  {10pt}   
  {10pt}   
  {\normalfont}  
  {}      
  {\bfseries} 
  {}     
  {.5em}  
  {}      
  
\theoremstyle{mylemma}
\newtheorem{lemma}{Lemma}[section]

\newtheoremstyle{numberfirst}
  {10pt}    
  {10pt}    
  {\normalfont} 
  {}       
  {\bfseries} 
  {:}      
  {.5em}   
  {\thmnumber{#2}\thmname{ #1}} 

\newtheorem{definition}{Definition}[section]

\usepackage{amsthm}

\theoremstyle{remark}
\newtheorem*{rem}{Remark}

\newcommand{\bE}{\mathbf{E}}

\newcommand{\bS}{\mathbf{S}}

\newcommand{\cP}{\mathcal{P}}
\newcommand{\cS}{\mathcal{S}}

\newcommand{\cC}{\mathcal{C}}

\begin{document}
\pagestyle{myheadings}
\markboth{\centerline{\'Aron Vil\'agi and Jen\H o Szirmai}}
{Power with Respect to Generalized Spheres Hyperbolic Geometry}

\title{Power with Respect to Generalized Spheres and Radical Surfaces in $\mathbf{H}^n$\thanks{Mathematics Subject Classification 2010: 51M10, 51M15, 52A20, 52C17, 52C22, 52B15. \newline
Key words and phrases: Hyperbolic geometry, power of a point, radical axis, hyperball packings, packing density.}}

\author{\'Aron Világi and Jen\H o Szirmai \\[1ex]
\normalsize Department of Algebra and Geometry, Institute of Mathematics,\\
\normalsize Budapest University of Technology and Economics, \\
\normalsize M\"uegyetem rkp. 3., H-1111 Budapest, Hungary \\
\normalsize szirmai@math.bme.hu}

\date{\normalsize{2026}}   


	\maketitle
\begin{abstract}
    \hyphenpenalty=10000
\noindent
This paper presents a unified theory for the power of a point with respect to generalized spheres (spheres, horospheres, and hyperspheres) in $n$-dimensional hyperbolic space $\mathbf{H}^n$. By extending the classical secant theorem, we derive a novel formula for hyperspheres and also prove that the radical surface of any two non-concentric generalized spheres is a hyperplane. These results provide tools for constructing power diagrams and studying hyperball \mbox{packings}.
\exhyphenpenalty=10000
\end{abstract}
\section{Introduction}

The concept of the power of a point with respect to a circle, originally introduced by Jakob Steiner in 1826 \cite{S1826}, is one of the cornerstones of classical Euclidean geometry.
It reflects the relative position of a given point with respect to a given circle.
While this concept is well-understood in Euclidean space $\mathbf{E}^n$ and has analogies in spherical geometry $\mathbf{S}^n$, its generalization to hyperbolic space $\mathbf{H}^n$ presents unique challenges and deep geometric structures.
In $n$-dimensional hyperbolic space, "the family of circles and spheres" is richer than in Euclidean geometry.
Besides the classical spheres, we must consider \textit{generalized spheres}: horospheres (surfaces orthogonal to a pencil of parallel lines) and hyperspheres (surfaces at a constant distance from a hyperplane).
While the geometry of classical hyperbolic spheres is widely discussed in the literature, the properties of hyperspheres~—~ \newline specifically regarding the power of a point and related incidence theorems ~—~ have received less attention, despite their growing importance in geometry.
The primary aim of this paper is to establish a unified theory for the power of a point with respect to generalized spheres in $\mathbf{H}^n$, with a special focus on the often-neglected case of hypercycles/hyperspheres.

\medskip
\textbf{Our main contributions are as follows:}

\begin{enumerate}
    \item Unlike spheres and horospheres, a hypersphere consists of two separate branches.
    We formulate and prove the \textit{Power of a point theorem for hypercycle ~-~ and hypersphere branches} (Theorem~\ref{thm:hypc}).
    This is a non-trivial extension of the classical theorem using secant lines intersecting both branches of a hypersphere, which gives a beautiful relation involving both the hyperbolic tangent ($\tanh$) and hyperbolic cotangent ($\coth$) functions.
    Using this theorem, we can define the power of any point with respect to a given hypersphere branch.
    \item We prove that the \textit{radical surface} of any two non-concentric generalized spheres (whether they are spheres, horospheres, or hyperspheres) is always a hyperbolic hyperplane in $\mathbf{H}^n$ (Lemma~\ref{lem:radax3}).
    This result is fundamental for constructing power diagrams (Dirichlet-Voronoi-like cells) in hyperbolic space (see \cite{Sz17, YSz25}).
\end{enumerate}

The motivation for this study also stems from the theory of ball packings and coverings (see \cite{Fe23}).
In recent years, determining the upper bound density of hypersphere packings in $\mathbf{H}^n$ has become a central topic of research \cite{Sz23,YSz25}.
To handle these complex arrangements, one often employs decomposition methods (e.g., dividing space into truncated simplices).
Such decompositions can be constructed using the radical surfaces of hyperspheres.

\medskip
A striking feature of our approach is its methodology.
Rather than relying on lengthy analytic computations or coordinate-heavy derivations, we employ \textit{synthetic geometric reasoning}.
By utilizing the properties of the Poincaré ball model and geometric inversion, we provide elegant, visual proofs that highlight the intrinsic beauty of hyperbolic geometry.

\medskip
\textbf{Structure of paper.} \; Section~2 and 3 review the necessary Euclidean preliminaries and the connection between the secant theorem and inversion.
Section~4 presents the main results regarding the power of a point in $\mathbf{H}^n$, including the new theorem on hypercycle branches.
Section~5 discusses the properties of hyperbolic inversion, and in Section~6, we prove that the radical surfaces are indeed hyperplanes.
Finally, we briefly discuss the implications of these results for non-congruent hyperball packings in $\mathbf{H}^n$.

\section{Preliminaries}

\subsection{The Euclidean case}

In the Euclidean plane, the power of a point $P$ with respect to a circle $\cC^E$ with center $O$ and radius $r$ is defined as
\[
\cP(P) := d^E(P,O)^2 - r^2,
\]
where $d^E(P,O)$ denotes the Euclidean distance between $P$ and $O$.
This definition was introduced by Jakob Steiner in \cite{S1826}, and it encodes the relative position of $P$ with respect to $\cC^E$.

\begin{thm}[Power of a point theorem in $\bE^2$]
Let a line through a point $P$ intersect the circle $\cC^E$ at points $A$ and $B$.
Then, independently of the choice of the line,
\[
d^E(P, A) \cdot d^E(P, B) = \mathcal{P}(P).
\]
If the line through $P$ is tangent to the circle $\cC^E$ at point $T$, then
\[
d^E(P, T)^2 = d^E(P, A) \cdot d^E(P, B).
\]
\end{thm}

The theorem admits natural analogues in higher dimensions: if $P$ is a point and $\cS^E$ is an $n$-dimensional sphere in $\bE^n$, then the power of $P$ with respect to $\cS^E$ can be expressed in the same way, and the product of secant segments through $P$ is constant.
\begin{rem} 
The power of a point theorem can be viewed as a unification of two classical results: when $P$ lies inside or on the circle, it becomes the theorem of intersecting chords, and when $P$ lies outside, it becomes the theorem of intersecting secants (both theorems can be found in the \textit{Elements} of Euclid).
\end{rem}
\begin{figure}[H]
	\centering
	\includegraphics[width=8.8cm]{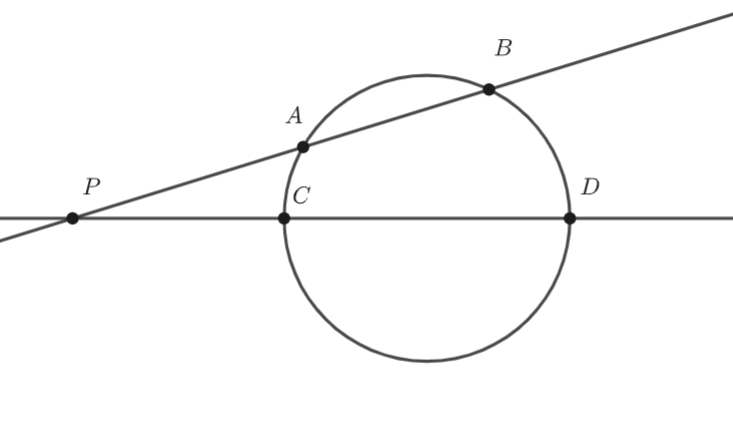}
	\caption{Figure~1: Power of a point theorem}
	\label{The secant-tangent theorem}
\end{figure}

\subsection{Relation with inversion}

The connection between the power of a point theorem and inversion provides a short proof:
Let $\cP(P) := d^E(P, C) \cdot d^E(P, D)$, where $CD$ is the diameter cut out from the circle $\cC^E$ with center $O$ by the line $PO$.
If $P$ lies outside $\cC^E$, inversion with respect to the circle centered at $P$ with radius $\sqrt{\cP(P)}$ leaves the Thales circle of the segment $CD$  (which is $\cC^E$) invariant, since inversion maps circles (not passing through the center of inversion) to circles.
Thus, it holds for every secant line that the product $d^E(P,A)\cdot d^E(P,B)$ is constant (and equals the power of point $P$ with respect to $\cC^E$), since, for example, the inverse image of $A$ is $B$ under this inversion.
\qed

If $P$ is located inside the circle $\cC^E$, a proof completely identical to the above can be given by applying negative inversion (in this case, a point and its image are located on opposite sides of the center of inversion).
The proof can be extended to higher dimensions using inversion with respect to an $n$-dimensional sphere.

\subsection{The spherical case}

An analogous statement holds in spherical geometry. Let $\cC^S$ be a circle in the unit sphere $\bS^2$ of radius $r$, and let $P$ be a point not antipodal to any point of $\cC^S$.

\newpage
\begin{thm}[Power of a point theorem in $\bS^2$]
Let a great circle through a point $P$ intersect the circle $\cC^S$ at points $A$ and $B$.
Then, independently of the choice of the line,
\[
\tan\!\left(\frac{d^S(P,A)}{2}\right)\cdot \tan\!\left(\frac{d^S(P,B)}{2}\right)
\]
If the great circle is tangent to $\cC^S$ at point $T$, then
\[
\tan^2\!\left(\frac{d^S(P,T)}{2}\right) = 
\tan\!\left(\frac{d^S(P,A)}{2}\right)\cdot 
\tan\!\left(\frac{d^S(P,B)}{2}\right).
\]
\end{thm}

This establishes the spherical analogue of the power of a point theorem.
The earliest works discussing this relationship on the spherical plane date back to the late-18th century (see \cite{LAJ}).
Since then, numerous new proofs of the theorem have been developed, and it is considered one of the fundamental relations in spherical trigonometry.

\section{Auxiliary lemmas}
   
In this section, we introduce two Euclidean lemmas that will serve as fundamental tools throughout the paper.

\begin{definition}
The \emph{radical axis} of two circles in the Euclidean plane is the locus of points having equal power with respect to both circles.
(This definition can be extended to $n$ dimensions for two $n$-dimensional spheres.)
\end{definition}

The name and the definition were given by Louis Gaultier 1813 in \cite{GL}.
It is well known that the radical axis of two non-concentric circles in the Euclidean plane is a line (and a hyperplane in higher dimensions).

\begin{lemma}\label{lem:important1}
Invert circles $k$ and $l$ with respect to a circle centered on the radical axis of $k$ and $l$.
Let $k'$ and $l'$ denote the inverse images.
Then the radical axis of $k$ and $l$ coincides with the radical axis of $k'$ and $l'$.
\end{lemma}

\begin{figure}[H]
	\centering
    \caption{Figure~2: Illustration of Lemma~\ref{lem:important1}}
	\includegraphics[width=7.1cm]{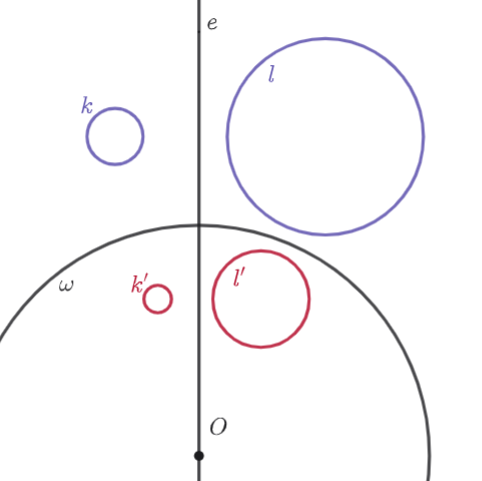}
    \label{pic:second}
\end{figure}

\paragraph{Proof:} Let $k$ and $l$ be circles, let $e$ denote their radical axis, and let $O$ be a point on line $e$.
Let $\omega$ be the circle of inversion with center $O$ and radius $r \in \mathbb{R}^+$.
By definition, $O$ lies on the radical axis of $k'$ and $l'$.

\medskip
Connect the centers of $k$ and $l$ with a line $f$.
This line is perpendicular to circles $k$ and $l$.
$f'$ (the inverse image of $f$ with respect to $\omega$) is an open circle passing through $O$ orthogonal to $k'$ and $l'$, with its center on $e$ (due to the properties of inversion).
Due to orthogonality, the center of $f'$ lies on the radical axis of $k'$ and $l'$.
Thus, the radical axis of $k'$ and $l'$ is the line connecting the center of $f'$ and $O$, which is the line $e$.
\qed

\begin{rem} 
The statement holds in $n$ dimensions for spheres and radical hyperplanes as well.
Since inversion is well-defined in $\mathbf{E}^n$ spaces and preserves all its important properties (see Chapter 2 of \cite{TWP}), the proof remains valid if $k$, $l$, $k'$, $l'$, $f'$ and $\omega$ are ($n$-dimensional) spheres, while $e$ and $f$ are hyperplanes.
Thus, we obtain the following lemma:
\end{rem}

\begin{lemma}\label{lem:important3}
Invert $n$-dimensional spheres $k$ and $l$ with respect to an $n$-sphere centered on the radical hyperplane of $k$ and $l$.
Let $k'$ and $l'$ be the inverse images.
Then the radical hyperplane of $k$ and $l$ coincides with the radical hyperplane of $k'$ and $l'$.
\end{lemma}
%
\begin{lemma}\label{lem:ortogcirc}
Let $\omega_1$ and $\omega_2$ be two circles intersecting orthogonally at points $A$ and $B$.
Let $c$ be a circle passing through $A$ and $B$.
Let $c'$ and $c''$ denote the inverse images of $c$ with respect to $\omega_1$ and $\omega_2$, respectively.
Then circles $c'$ and $c''$ coincide.
\end{lemma}
\paragraph*{Proof:} The lemma follows directly from the conformality of inversion.
\smallskip

\begin{figure}[H]
	\centering
	\caption{Figure~3: Illustration of Lemma~\ref{lem:ortogcirc}}
	\includegraphics[width=10.2cm]{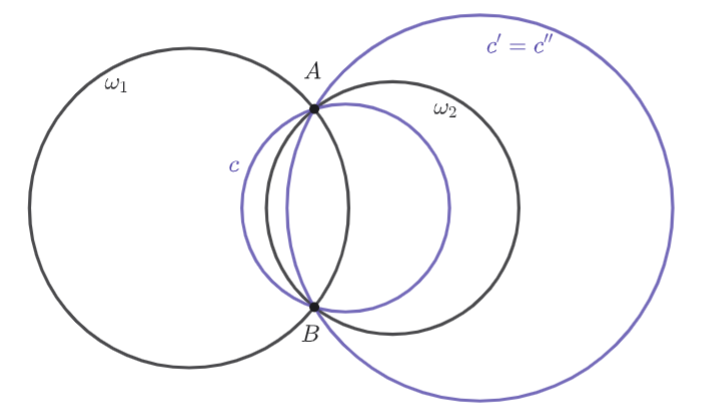}
	\label{pic:third}
\end{figure}

\section{The notion of power hyperbolic geometry}

\subsection{The power of a point theorem in hyperbolic geometry}

\begin{thm}[Power of a point theorem for circles in $\mathbf{H}^2$]\label{thm:hppt}
Let a line through a given point $P$ intersect the circle $\cC^H$ at $A$ and $B$ in the hyperbolic plane.
Then the following product is constant for every line that intersects $\cC^H$:
	\[\mathrm{tanh}\left(\frac{d^H(P,A)}{2}\right) \cdot \mathrm{tanh}\left(\frac{d^H(P,B)}{2}\right).
\]
    \end{thm}

\begin{figure}[H]
	\centering
	\caption{Figure~4: Hyperbolic secant theorem in the Poincaré disk model}
	\includegraphics[width=6cm]{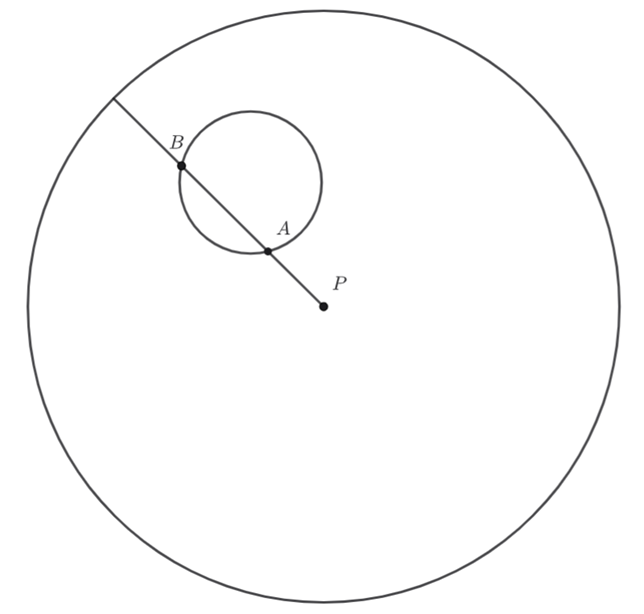}
	\label{pic:fourth}
\end{figure}

\paragraph*{Proof:} 
The proof uses the Poincaré disk model of the hyperbolic plane. It is sufficient to examine the case where $P$ lies at the center of the model and $\mathcal{C}^H$ is arbitrary, since any other case can be mapped to this one by a suitable isometry. Placing point $P$ at the center of the model is particularly useful because (Euclidean) lines passing through the center in the Poincaré disk model correspond to lines of the hyperbolic plane. Furthermore, if the Euclidean distance of a point from the center is $r$ in the model, its hyperbolic distance in the hyperbolic plane is $2\,\mathrm{artanh}\,r$.

\medskip
When $P$ lies at the center of the model, it follows that the Euclidean lengths of the secant segments are $d^E(P, A) = \mathrm{tanh}\!\left(\frac{d^H(P, A)}{2}\right)$ and $d^E(P, B) = \mathrm{tanh}\!\left(\frac{d^H(P, B)}{2}\right)$.

\medskip
Since the secant theorem holds in the Euclidean plane, and circles of the hyperbolic plane are represented by circles in the model, the product
\[
\mathrm{tanh}\!\left(\frac{d^H(P, A)}{2}\right)\cdot \mathrm{tanh}\!\left(\frac{d^H(P, B)}{2}\right)
\]
is constant for every line intersecting $\mathcal{C}^H$. \hfill $\square$

\medskip
If $P$ lies \emph{outside} the circle $\cC^H$, it follows similarly that if the line passing through $P$ is tangent to $\cC^H$ at point $T$, then
\[
\mathrm{tanh}\!\left(\frac{d^H(P, A)}{2}\right) \cdot \mathrm{tanh}\!\left(\frac{d^H(P, B)}{2}\right) 
= \mathrm{tanh}^2\!\left(\frac{d^H(P, T)}{2}\right).
\]

\medskip
The same theorem applies to horocycles. Since a horocycle appears in the Poincaré disk as a circle tangent to the boundary circle, the proof is identical to the one presented above.
The relationship remains valid for secant segments intersecting the horocycle "at infinity" if we consider $\mathrm{tanh}(\infty)$ to be $1$.

\begin{rem}
These statements hold for $n$-dimensional spheres and horospheres as well.
The proof is identical to the planar case if we use the $n$-dimensional Poincaré ball model (for the $n$-dimensional model, see \cite{BCSS}), since the secant theorem is valid in the $n$-dimensional Euclidean space.
\end{rem}

\begin{definition}
The \emph{power} of a point $P$ with respect to a circle or horocycle in the hyperbolic plane is defined by the product:
\[
\mathrm{tanh}\!\left(\frac{d^H(P, A)}{2}\right)\cdot \mathrm{tanh}\!\left(\frac{d^H(P, B)}{2}\right).
\]
The power with respect to an $n$-dimensional sphere or horosphere is defined in the same way.
\end{definition}

\subsection{Power with respect to hypercycles and hyperspheres}
Since hypercycles in the Poincaré disk model are circular arcs (or line segments) that are not perpendicular to the boundary circle, it is easy to see that the hyperbolic power of a point theorem can be extended to a branch of hypercycle, if we can intersect that branch twice with the secant line.
(This can be proved similarly to Theorem~\ref{thm:hppt}). From a given point $P$, however, a line can intersect a branch of a hypercycle twice only if that branch lies closer to $P$ than the other.\footnote{except when $P$ lies on the base line.
In that case both branches can be intersected twice at their points at infinity, if the secant is the base line itself}. Moreover, for any point and any hypercycle it is possible to draw a secant that intersects the hypercycle on two distinct branches.
In this case, the hyperbolic power of a point theorem no longer holds for the lengths of the secant segments.
To address these difficulties, we introduce a version of the hyperbolic power of a point theorem generalized for hypercycle branches.
In order to construct the theorem, we shall examine an interesting property of the Poincaré disk model:
\begin{lemma}\label{lem:hypc1}
Let $H_1$ and $H_2$ be two branches of a hypercycle in the hyperbolic plane.
In the Poincaré disk model, the inverse of the circular arc representing $H_1$ with respect to the boundary circle completes the circular arc representing $H_2$, and the inverse of the arc representing $H_2$ with respect to the boundary circle completes the arc representing $H_1$.
\end{lemma} 
\medskip

\begin{figure}[H]
	\centering
	\caption{Figure~5: Illustration of Lemma~\ref{lem:hypc1}}
	\includegraphics[width=10cm]{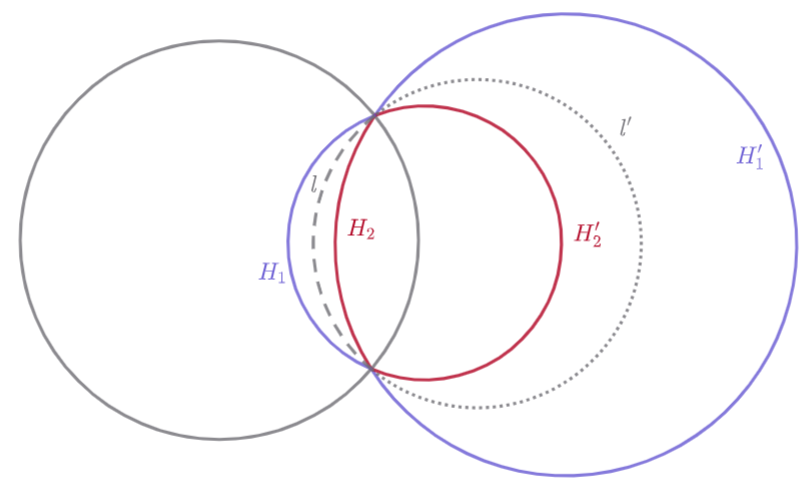}
	\label{pic:fifth}
\end{figure}
\medskip
\paragraph{Proof:} 
Since inversion with respect to a circular arc representing a line in the Poincaré disk model models reflection, by the definition of a hypercycle, the inverse of the arc representing $H_1$ with respect to the arc of the base line $l$ is the arc representing $H_2$.
Furthermore, the arc of the base line $l$ in the model is orthogonal to the boundary circle of the model;
thus, by Lemma~\ref{lem:ortogcirc}, the full circles of the arcs modeling $H_1$ and $H_2$ are inverses of each other with respect to the boundary circle.
Since the arcs corresponding to $H_1$ and $H_2$ lie inside the boundary circle, their inverse images with respect to the boundary circle lie outside the model, completing each other's arcs into full circles. \qed

\medskip
The lemma naturally extends to the $n$-dimensional Poincaré ball model, where $H_1$ and $H_2$ are $n$-dimensional hyperspheres, as they appear as $n$-dimensional spherical caps in the model.
The proof proceeds analogously to the one above, since Lemma~\ref{lem:ortogcirc} is also true for spheres of higher dimensions.

\begin{thm}[Power of a point theorem with respect to hypercycle-branches]\label{thm:hypc}
Let $H$ be a hypercycle in the hyperbolic plane with branches $H_1$ and $H_2$.
Let $P$ be a point in the hyperbolic plane.
\begin{enumerate}
    \item[\textup{(i)}] \textbf{Secant on the same branch:}  
    If a hyperbolic line passing through $P$ intersects branch $H_1$ at $A$ and $B$ (provided this is possible), then the product
    \[
    \tanh\!\left(\frac{d^H(P, A)}{2}\right)\cdot \tanh\!\left(\frac{d^H(P, B)}{2}\right)
    \]
    is constant, independently of the choice of the secant line passing through $P$.
    \item[\textup{(ii)}] \textbf{Secant on different branches:}  
    If a hyperbolic line passing through $P$ intersects branch $H_1$ at $A$ and the opposite branch $H_2$ at $B$, then the product
    \[
    \tanh\!\left(\frac{d^H(P, A)}{2}\right)\cdot \coth\!\left(\frac{d^H(P, B)}{2}\right)
    \]
    is constant, independently of the choice of the secant line.
    Moreover, this product equals the product in case (i) if $H_1$ can be intersected at two points by a secant line passing through $P$.
\end{enumerate}
\end{thm}

 \begin{rem}
 Equivalently, since $\coth\left(\frac{PB}{2}\right) = \tanh\left(\frac{PB + i\pi}{2}\right)$, to keep $\tanh$ in the formula, one can also write:
    \[
        \tanh\left(\frac{d^H(PA)}{2}\right) \cdot \tanh\left(\frac{d^H(PB) + i\pi}{2}\right).
    \]
\end{rem}

\paragraph{Proof:}
Place point $P$ at the center of the Poincaré disk model and consider the full circle of the arc representing $H_1$.
Let the line drawn from $P$ through $A$ intersect this circle at another point $A_2$ in the Euclidean plane of the model.
By Lemma~\ref{lem:hypc1}, the arc completing the arc modeling $H_1$ outside the model is the inverse of the arc modeling $H_2$ with respect to the boundary circle (which is a unit circle).
Since $\coth{x}$ and $\tanh{x}$ are reciprocals, $\coth(\frac{d^H(P, B)}{2})$ corresponds to the distance $d^E(P, A_2)$ in the Euclidean plane of the model if $B$ lies on $H_2$ (see Figure 6).
Otherwise, $A_2$ coincides with $B$.
Since the power of a point theorem holds in the Euclidean plane, the statement is proven.\qed

\medskip
The theorem applies equally to branches of $n$-dimensional hyperspheres. The proof is identical to the one above, as Lemma~\ref{lem:hypc1} holds in higher dimensions as well.

\begin{figure}[H]
	\centering
	\caption{Figure~6: Illustration of Proof for Theorem~\ref{thm:hypc}}
	\includegraphics[width=9cm]{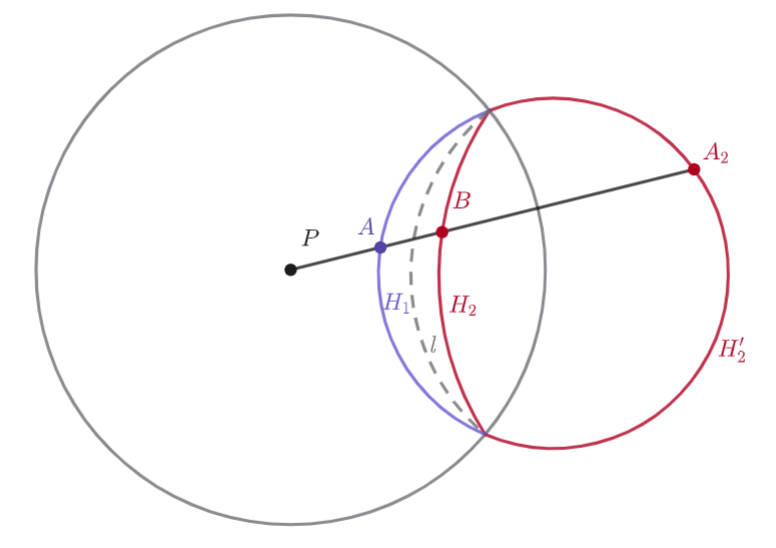}
	\label{pic:sixth}
\end{figure}
Based on the theorem above, we can provide a more general definition for the concept of power with respect to hypercycle branches:

\begin{definition}
    Let $P$ be a point and $H$ a hypercycle in the hyperbolic plane with branches $H_1$ and $H_2$.
    Draw an arbitrary secant line through $P$ intersecting branches $H_1$ and $H_2$ at points $A$ and $B$, respectively (such a secant always exists).
    The \emph{power of point $P$ with respect to $H_1$} is determined by the following product:
    \[
    \tanh\!\left(\frac{d^H(P, A)}{2}\right)\cdot \coth\!\left(\frac{d^H(P, B)}{2}\right).
\]
    The power with respect to branches of higher-dimensional hyperspheres is defined by the same expression.
\end{definition}

\begin{rem} By this definition, the powers of a given point $P$ with respect to the two branches of a given hypercycle are reciprocals of each other.
\end{rem}
\begin{rem}
Based on the above, it is easy to see that if $P$ lies on the base line of the hypercycle, its power with respect to both branches equals $1$.
If $P$ lies on the branch for which the power is being examined, the power is $0$.
If $P$ lies on the other branch, the power can be considered $\infty$, since the right-hand limit of the $\coth$ function at $0$ is infinity.
Following the proof, this also has a nice geometric representation.
If the arc representing $H_2$ passes through the center of the Poincaré disk, the arc representing $H_1$ appears as a straight line segment.
\end{rem}

\section{Inversion in spaces of constant curvature}
    \begin{definition}[Inversion in spherical geometry]
    The inverse of a point $P \in \mathbf{S}^2$ with respect to a circle $\cC^S$ with center $O$ and radius $r$ is the point $P' \in \mathbf{S}^2$ such that:
    $$\tan\left({\frac{d^S(O,P)}{2}}\right) \cdot \tan\left({\frac{d^S(O, P')}{2}}\right) = \tan^2\left({\frac{r}{2}}\right)$$
    (where  $d^S(O,P')$ is measured in the direction of $P$).
\end{definition}
    This definition is similar to the definition of inversion in the Euclidean plane, where the defining equation is $d^E(O,P) \cdot d^E(O,P') = r^2$.
    The Euclidean and the spherical inversion share all of their most important properties (see \cite{CJ}).
    We will discuss the inversion in the hyperbolic plane in more detail, since this case is less well-known than the Euclidean and spherical cases.
    \begin{definition}[Inversion in the hyperbolic plane]
    The inverse of a point $P \in \mathbf{H}^2$ with respect to a circle  $\cC^H$ with center $O$ and radius $r$ is the point $P' \in \mathbf{H}^2$ such that:
    $$\tanh\left({\frac{d^H(O,P)}{2}}\right) \cdot \tanh\left({\frac{d^H(O, P')}{2}}\right) = \tanh^2\left({\frac{r}{2}}\right)$$
    (where $P'$ is lying on the ray from $O$ through $P$).
\end{definition}
	This transformation has a domain of definition, since the range of $\tanh$ on positive values is not all positive numbers, 
	only the interval $(0,1)$.
    Hence, points within a distance less than or equal to $2\mathrm{artanh}\left(\tanh^2\left(\frac{r}{2}\right)\right)$ 
	from $O$ will not have a real image under this transformation.
    So the transformation is undefined on the closed disk of radius 
	$2\mathrm{artanh}\left(\tanh^2\left(\frac{r}{2}\right)\right)$ centered at $O$.
    (This can be verified algebraically.)
	
	\subsubsection*{5.3 \;\;Properties of the transformation}
    The properties of the inversion in the hyperbolic plane is very similar to the properties of the Euclidean inversion:
	\begin{enumerate}[(1.)]
		\item The domain and codomain are equal, and the mapping is bijective between them.
		\item A circle within the domain is mapped to another circle.
		\item A line is mapped to an arc orthogonal to the boundary circle of the domain of definition.
		(and vice versa).
		\item A hypercycle is mapped to an arc (non-orthogonally) intersecting the boundary circle of the domain of definition.
		(and vice versa).
		\item Angle-preserving.
	\end{enumerate}
	\paragraph{Proofs:}
	\begin{enumerate}[(1.)]
		\item Follows from the commutativity of multiplication.
		\item Use the Poincaré disk model.
		Again, it is enough to examine the case where $O$ is the center of the model.
		When $O$ is the center point, the transformation can be modeled as inversion in a Euclidean circle centered at $O$ with radius $\tanh\left(\frac{r}{2}\right)$ (this can be verified algebraically).
		Since circles not passing through the center are mapped to circles under inversion, and since Euclidean circles in the model correspond to hyperbolic circles (and vice versa), the statement follows.
		\item Follows similarly to (2.) from the properties of Euclidean inversion.
		\item Follows similarly to (2.) from the properties of Euclidean inversion.
		\item Follows from the conformality of the Poincaré model and Euclidean inversion.
	\end{enumerate}
    
    \begin{rem} One key property of classical inversion does not hold here: there do not exist $O$ and $r$ for which any two  (non-congruent) circles are mapped into one another. However, this statement is true for any two circles whose external common tangents intersect (this can be easily verified).
    \end{rem}
    In Euclidean geometry, as mentioned earlier, inversion with respect to an $n$-dimen- sional circle is defined the same way as in the planar case, and has the $n$-dimensional equivalents of the properties of the planar version.
    In the $n$-dimensional Poincaré ball model of hyperbolic geometry, the transformation representing the reflection of a hyperplane is a Euclidean inversion with respect to an $n$-dimensional sphere (which represents the hyperplane in the ball model).

    \medskip
    We can define the inversion in $n$-dimensional hyperbolic spaces with respect to an $n$-dimensional sphere the same way as in the planar case preseving all important properties.
    The proofs for these are the same as those given for the planar case, only using the $n$-dimensional inversion and the $n$-dimensional Poincaré ball model.

    \begin{figure}[H]
	\centering
	\caption{Figure~7: Hyperbolic inversion with respect to $\omega$ in the Poincaré disk model}
	\includegraphics[width=7.81cm]{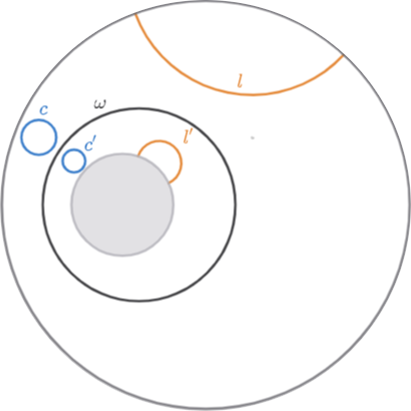}
	\label{pic:seventh}
\end{figure}

\section{ Radical axes and surfaces in hyperbolic geometry}
    
	As shown previously, the hyperbolic and spherical analogues of the power of a point theorem allows us to define the power of a point with respect to a circle in these geometries.
    (This concept can be extended to generalized circles as well in the hyperbolic case.) Thus, the definition of the radical axis in these geometries coincides with the definition given in Euclidean geometry.
    In the spherical case, it is well known that the radical axis of two circles is an arc of a great circle, in other words, a spherical line.
    (See in \cite{ARE}.) In the following we examine the radical axis of circles, hypercycles and horocycles on the hyperbolic plane and also provide a generalization to $\mathbf{H}^n$ spaces.

    \subsection{Radical axis of two circles in hyperbolic geometry}
    \begin{lemma}\label{lem:radax1}
    The radical axis of any two circles in hyperbolic geometry is a hyperbolic line.
    (Or an empty set of points in concentric cases.)
    \end{lemma}
	\paragraph{Proof:}
	Let $k$ and $l$ be any two circles in the hyperbolic plane.
    Naturally, there exist (infinitely many) points having equal power with respect to both.
    Take such a point $O$, and use the Poincaré disk model where $O$ is at the center.
    Since all lines through the origin are also straight lines in the model, $O$ lies (in the Euclidean sense) on the radical axis of the Euclidean images of $k$ and $l$.
    The Euclidean radical axis of these circles, denoted $e$, is then also a hyperbolic straight line.

    \medskip
	Choose another point $P \ne O$ on this line. Using an isometry (e.g. reflection), map $P$ to the center of the model.
    This transformation in the model is an inversion in a circle centered on $e$ and orthogonal to the boundary circle.
    By Lemma~\ref{lem:important1}, the image circles $k'$, $l'$ will also have radical axis $e$.
    Since $P$ now lies at the origin, by the reasoning used in the proof of Theorem~\ref{thm:hppt}, $P$ must also lie on the hyperbolic radical axis of $k'$ and $l'$.
    Therefore, every point on $e$ has equal power with respect to $k$ and $l$. By continuity, no other points do.
    Hence, the radical axis is this line $e$.

    \medskip
	In the concentric case the radical axis is trivially an empty set similarly to euclidean geometry.
    \qed

        \begin{figure}[H]
		\centering
		\caption{Figure~8: Radical axis of two hyperbolic circles in the Poincaré disk model}
		\includegraphics[width=6.9cm]{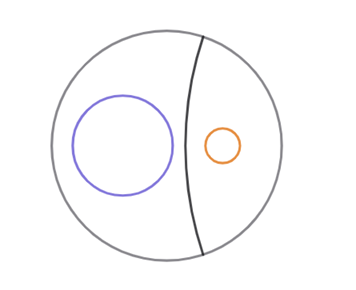}
		\label{pic:sixth}
	\end{figure}

	\begin{rem}
	If we take the intersection of the radical axis and the hyperbolic line connecting the centers of the two circles and place it at the center of the Poincaré model, we can see by symmetry that the radical axis is perpendicular to this line.
    \end{rem}

    \subsection{Generalization to cycles:}
	By extending the reasoning used in the proof of Theorem~\ref{thm:hppt} to arbitrary cycles, we can conclude that a point $P$ lies on the radical axis of any two cycles (circle, horocycle, or hypercycle branch) if and only if, when placing point $P$ at the center of the Poincaré disk model, $P$ lies on the Euclidean radical axis of the circles modeling the two cycles (circles or circular arcs).
    Thus, in a completely analogous way to the proof of Lemma~\ref{lem:radax1}, we obtain the following statement:
    \begin{lemma}\label{lem:radax2}
    The radical axis of two non-concentric generalized circles in the hyperbolic plane is also a hyperbolic line.
    \end{lemma}
	
	\begin{rem} It can be easily seen by the methods used above that in cases where two hypercycles have the same baseline, the radical axis is the baseline itself, and that where two horocycles have the same infinite point, the radical axis of these cycles is an empty set of points.
    \end{rem}
    \begin{rem} The three radical axes determined by three circles/cycles in the hyperbolic plane meet is one point.
    This can be proved analougusly to the Euclidean case.
    \end{rem}

    \begin{figure}[H]
		\centering
		\caption{Figure~9: Radical axes of a circle, a hypercycle and a horocycle}
		\includegraphics[width=8cm]{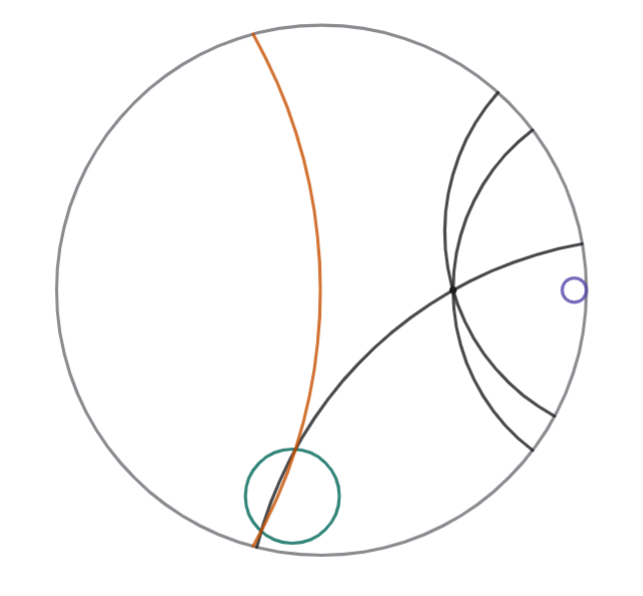}
		\label{pic:sixth}
	\end{figure}
	
	\subsection{Generalization to $n$ dimensions:}
	All the statements used in the above proofs extend naturally to higher dimensions.
    The proofs are identical to the proofs for Lemma~\ref{lem:radax1} and Lemma~\ref{lem:radax2} using the \ref{lem:important3} version of Lemma~\ref{lem:important1} and the $n$-dimensional Poincaré ball-model of the $\mathbf{H}^n$ spaces.
    So we get:
    \begin{lemma}\label{lem:radax3}
    The radical surface of two (non-concentric) generalized spheres in the $n$-dimensional hyperbolic space is a hyperplane.
    \end{lemma}
\medskip
\section{Future perspectives on non-congruent hyperball packings}

In $n$-dimensional hyperbolic geometry, numerous new questions arise regarding packing and covering problems.

\medskip
In the space $X^n$, let $d_n(r) \;~(n\geq 2)$ denote the density of $n+1$ mutually tangent spheres or horospheres (in the case of $r=\infty$) of radius $r$ relative to the simplex spanned by their centers.
L. Fejes Tóth and H. S. M. Coxeter conjectured that the packing density of spheres of radius $r$ in $X^n$ cannot exceed $d_n(r)$.
The conjecture was proved by C. A. Rogers for the Euclidean space $\mathbf{E}^n$, and the two-dimensional spherical case was solved by L. Fejes Tóth.

\medskip
The maximum density in $\mathbf{H}^3$ is $\approx 0.85328$, which is realized by horosphere packings, 
the regular tetrahedral tessellation corresponding to the ideal packing is given by the Coxeter-Schläfli symbol $\{3,3,6\}$.
Sphere packings in hyperbolic $n$-space and other Thurston geometries are widely discussed in the literature.
However, while there are countless open questions currently being investigated regarding sphere and horosphere packings and coverings, relatively few results are known concerning hyperball arrangements.

\medskip
In the hyperbolic plane $\mathbf{H}^2$, the universal upper bound for hypercycle packing density is $\frac{3}{\pi}$, proved by I. Vermes.
Similarly for hypercycle coverings, the universal lower bound for covering density is $\frac{\sqrt{12}}{\pi}$.

\medskip
The following papers deal with higher-dimensional congruent hyperball packings and coverings.

\medskip
\cite{Sz17} examines congruent hyperball packings in $3$-dimensional hyperbolic space and presents a decomposition algorithm that ensures the decomposition of $\mathbf{H}^n$ into truncated tetrahedra for every saturated hyperball packing.
Therefore, to obtain the upper bound for the density of hyperball packings, it is sufficient to determine the upper bound for the density of congruent hyperball packings within truncated simplices, which is $\approx 0.86338$.

\medskip
\cite{YSz25} continues the investigation of congruent hyperball packings in higher-dimen- sional hyperbolic spaces $\mathbf{H}^n$ ($n\ge4$) and shows that for every $n$-dimensional congruent, saturated hyperball packing, there exists a decomposition of the $n$-dimensional hyperbolic space into truncated simplices.
Furthermore, the paper proves that the upper bound for the density of saturated congruent hyperball packings associated with the corresponding truncated tetrahedron cells is realized in a regular truncated tetrahedron.
In $4$-dimensional hyperbolic space, this density upper bound is $\approx 0.75864$ (which is larger than the previously conjectured maximum).
Moreover, it refutes the conjecture of A. Przeworski regarding the monotonicity of the congruent hyperball packing density function in $4$-dimensional hyperbolic space.

\medskip
In $\mathbf{H}^n$ space, cell decompositions defined by the radical hyperplanes of hyperballs provide an opportunity to develop methods similar to the above and to determine the density upper bound for incongruent hyperball packings.

\end{document}